\documentclass[12pt]{article} 
\usepackage{latexsym}
\usepackage{amssymb}
\usepackage{euscript}
\usepackage{color} 

\let\cal\mathcal

\usepackage{latexsym}
\usepackage{amssymb}
\usepackage{euscript}
\usepackage[dvips]{graphics}
\usepackage{epsf}

\let\cal=\mathcal      

\def\mcc{M\raise.5ex\hbox{c}C}
\def\mccarthy{M\raise.5ex\hbox{c}Carthy}

\def\l{\lambda}

\def\vare{\varepsilon}

\def\={\ = \ }

\def\E{E_\l}

\def\C{\mathbb C}
\def\R{\mathbb R}
\def\T{\mathbb T}
\def\D{\mathbb D}

\def\dis{\displaystyle}

\def\be{\setcounter{equation}{\value{theorem}} \begin{equation}}
\def\ee{\end{equation} \addtocounter{theorem}{1}}
\def\beq{\begin{eqnarray*}}
\def\eeq{\end{eqnarray*}}

\def\vs{\vskip 5pt}
\def\bs{\vskip 12pt}

\def\bp{{\sc Proof: }}
\def\ep{{}{\hfill $\Box$} \vskip 5pt \par}

\def\bl{\begin{lemma}}
\def\el{\end{lemma}}
\def\bt{\begin{theorem}}
\def\et{\end{theorem}}
\def\bprop{\begin{prop}}
\def\eprop{\end{prop}}
\def\bd{\begin{definition}}
\def\ed{\end{definition}}
\def\br{\begin{remark}}
\def\er{\end{remark}}
\def\bexer{\begin{exercise}}
\def\eexer{\end{exercise}}

\newtheorem{theorem}{Theorem}[section]
\newtheorem{prop}[theorem]{Proposition}
\newtheorem{lemma}[theorem]{Lemma}
\newtheorem{cor}[theorem]{Corollary}

\newtheorem{definition}[theorem]{Definition}

\def\Flip#1{1/\overline{#1}} 

\def\empty{\emptyset}
\def\n{\noindent}
\def\tnsym{$\T^n$-symmetric}
\def\cstarn{(\C^*)^n}

\newcommand{\U}{\cal U}

\def\E{\mathbb E}

\newcommand{\hib}{H^\infty (\D^2)}
\newcommand{\cn}{\C^n}
\newcommand{\cnm}{\C^{n-1}}
\newcommand{\dn}{\D^n}

\newcommand{\en}{\E^n}

\newcommand{\tn}{\T^n}
\newcommand{\tnm}{\T^{n-1}}

\newcommand{\inn}{\,\in\,}
\newcommand{\czn}{\C [ z_1, \dots, z_n]}
\newcommand{\cznm}{\C [ z_1, \dots, z_{n-1}]}
\newcommand{\va}{algebraic set}

\newcommand{\tp}{{p^{\sim}}}
\newcommand{\tq}{{q^{\sim}}}

\begin{document}
\setlength{\baselineskip}{21pt}
\title{Toral Algebraic Sets and Function Theory on Polydisks}
\author{Jim Agler
\thanks{Partially supported by National Science Foundation Grant
DMS 0400826}\\
U.C. San Diego\\
La Jolla, CA 92093
\and
John E. M\raise.5ex\hbox{c}Carthy
\thanks{Partially supported by National Science Foundation Grant
DMS 0501079}\\
Washington University\\
St. Louis, MO 63130
\and
Mark Stankus\\
California Polytechnic State University\\
San Luis Obispo, CA 93407}
\date{\today }

\bibliographystyle{plain}

\maketitle
\begin{abstract}
A toral algebraic set $A$ is an algebraic set in $\C^n$ whose
intersection with $\T^n$ is sufficiently large to determine the 
holomorphic functions on $A$. We develop the theory of these
sets, and give a number of applications to function theory in
several variables and operator theoretic model theory. 
In particular, we show that the uniqueness set for an extremal Pick problem
on the bidisk is a toral algebraic set, that rational inner
functions have zero sets whose irreducible components are not toral,
and that the model theory for a commuting pair of contractions with
finite defect lives naturally on a toral algebraic set.
\end{abstract}

%
\def\new#1{{#1}}
%

\baselineskip = 18pt

\setcounter{section}{-1}
\section{Introduction}

Throughout this paper, we shall let $\D$ denote the unit disk in
the complex plane, $\T$ be a unit circle, and let $A(\D^n)$
denote the polydisk algebra, 
the algebra of functions that are continuous on the closure of
$\D^n$ and holomorphic on the interior.

When studying function theory on the polydisk $\D^n$, it is often
useful to focus on the torus $\T^n$, which is the distinguished
boundary of $\D^n$. In several important ways, the behavior of a
function in $A(\D^n)$ is controlled by its behavior on $\T^n$: not
only is $\T^n$ a set of uniqueness,  but every function in the algebra
attains its maximum modulus on $\T^n$. 

Consider now some algebraic set $A$ contained in $\C^n$. When
studying function theory on $A$, the intersection of the torus with
$A$ may or may not play an important role. For example, if $A =
\{ (z,\dots, z ) \ : \ z \, \in \, \C \}$, then $A$ is a plane, 
and $A \cap \T^n$ is a unit circle. However, if 
$A = \{ (z,0,\dots, 0 ) \ : \ z \, \in \, \C \}$, then 
$A \cap \T^n$ is empty.

We shall say that a variety (by which we always mean an irreducible
algebraic set) $V$ is {\em toral} if its intersection
with $\T^n$ is fat enough to be a determining set for holomorphic
functions on $V$ (see Section~\ref{secc} for a precise
definition). Otherwise we shall 
call the variety {\em atoral}. We shall say that a polynomial in
$\czn$ is toral (respectively, atoral) if the zero set of every
irreducible factor is toral (respectively, atoral).

It turns out that factoring polynomials into their toral and atoral
factors is extremely useful when studying function theory on
$\D^n$.

Consider first the Pick problem on the bidisk, $\D^2$.
Let $\hib$ denote the Banach algebra of bounded analytic functions
on the bidisk.
A {\em solvable Pick problem on $\D^2$} is a set $\{\l_1, \dots,
\l_N \}$ of points in $\D^2$ and a set $\{w_1, \dots, w_N\}$ of
complex numbers such that there is some function $\phi$ of norm
less
than or equal to one
in $\hib$ that interpolates (satisfies $\phi(\l_i) = w_i \ \forall
\ 1 \leq i \leq N$).
An {\em extremal Pick problem} is a solvable Pick problem for which
no function of norm less than one interpolates.

Unlike extremal Pick problems on $\D$, an extremal Pick problem on $\D^2$
need not have a unique solution. However there is some subset 
$\U$ of $\D^2$, called the uniqueness set of the problem,
on which all solutions must agree ($\U$ obviously
contains all of the points $\l_i$). We prove in Theorem~\ref{thmf2}
that the 
uniqueness set
equals the intersection of $\D^2$ and the zero set of a toral polynomial. 

A second place where the concept of toral/atoral polynomials arises
is in the study of rational inner functions on $\D^n$.
An inner function is a function $\phi$ that is holomorphic and
bounded on $\D^n$ and whose radial boundary values, which exist
almost everywhere \cite[Thm. 3.3.5]{rud69}, have modulus one almost
everywhere. 
W.~Rudin showed \cite[Thm. 5.2.5]{rud69}
that every rational inner function can be represented in the form
\be
\label{eqa1}
\phi(z)
\= z^h \frac{\overline{p(1/\overline{z})}}{p(z)}
\ee
for some polynomial $p$ that does not vanish on $\dn$, and
some monomial $z^h$ such that $h \ge \deg p$.
We show in Theorem~\ref{thmd1} that in the representation (\ref{eqa1}),
the atoral factor of $p$ is uniquely determined, and the toral
factor is completely arbitrary.
As a consequence of this analysis, we show in Proposition~\ref{propd2}
that the zero set of a rational inner function is an atoral
algebraic set. In Proposition~\ref{propd3} we show that the
singular set of a rational inner function, namely the set of points
on $\T^n$ to which the function cannot be continuously extended
from $\D^n$, is always of dimension at most $n-2$.

The Sz.Nagy-Foia\c{s} model theory for a pair of commuting
contractions $(T_1,T_2)$ realizes them as the compression of a pair of
commuting isometries $(S_1,S_2)$ \cite{szn-foi}. 
In the event that one of the
isometries has finite defect, it can be represented as
multiplication by the independent variable on a vector-valued Hardy
space, and the other isometry becomes multiplication by a
matrix-valued inner function $B$. 
The set $$
A \= \{ (z,w) \, \in \, \C^2 \ : \ \det(B(z) - wI) = 0 \}
$$
is toral, and the extension spectrum of $(S_1,S_2)$ is $A \cap
\T^2$. This means that in addition to the $\D^2$ functional
calculus, one has an $A$ functional calculus, which is stronger.
In other words, for every polynomial $p$, instead of And\^o's inequality
\cite{and63}
$$
\| p(T_1, T_2) \| \ \leq \ \|p \|_{\D^2} ,
$$
one has
$$
\| p(T_1, T_2) \| \ \leq \ \|p \|_{A} .
$$
Therefore the study of the function theory of toral algebraic sets
is important in understanding the functional calculus for commuting
pairs of matrices. See the papers \cite{bv00} and \cite{bv03} 
by J.~Ball and
V.~Vinnikov and \cite{bsv04} by J.~Ball, C.~Sadosky and
V.~Vinnikov for another viewpoint.

The lay-out of this paper is as follows. In Section~\ref{secc}
we give precise definitions of the concepts of toral and atoral,
and make some basic observations. In Section~\ref{secc2}, we study
how various geometric properties are related to the analytic notion
of torality. In Section~\ref{secd} we give applications to the
study of rational inner functions, and in Section~\ref{secf} we
characterize the uniqueness set for a Pick problem on the bidisk.

We have tried to write our paper to be intelligible to both
algebraic geometers and to analysts. Consequently we apologize for
belaboring points that will be obvious to some readers.

\section{Toral and Atoral Algebraic sets in $\C^n$}
\label{secc}

For $p \inn
\czn$,  the polynomials in the commuting variables
$z_1,\ldots,z_n$ with complex coefficients,
we shall denote the zero set of $p$
by $Z_p$.
An {\em \va} in 
$\cn$ is a finite intersection of zero sets of
polynomials. If the algebraic set cannot be written as a 
finite union of strictly
smaller algebraic sets, it is called a {\em variety}.
Every algebraic set
$A$ can be written as a union of
varieties, no one containing another,
called  the
irreducible components of $A$.

The algebra of holomorphic functions on an algebraic set $A$, denoted
$Hol(A)$, consists of all complex-valued functions $f$ on $A$ with
the property that, for every point $\l$ in $A$, there is an open
set $U$ in $\cn$ containing $\l$, and a holomorphic function $g$
on $U$, such that $ g |_{U \cap A} = f |_{U \cap A}$.

Let $X \subseteq \cn$, and let $A$ be an algebraic set in $\cn$. 
We shall say that
{\em $X$ is determining for $A$} if $f \equiv 0$ 
whenever $f \inn Hol(A)$ and
$f |_{X \cap A} = 0$. 
Note that $X$ is not required to be a subset of $A$, which is a
departure from the usual use of ``determining''. However, it is
also determining in the following sense.

\bprop
Let $X \subseteq \cn$, and let
$A_1$ and $A_2$ be algebraic sets for which $X$ is determining.
If
$A_1 \cap X = A_2 \cap X$, then $A_1=A_2$.
\eprop

\bp
There exist nonzero polynomials $p_1,\ldots,p_m$ such that
$A_1 = Z_{p_1} \cap Z_{p_2} \cap \cdots \cap Z_{p_m}$.
For each $j$, $p_j$ vanishes on $A_1 \cap X$ (since
$p_j$ vanishes on $Z_{p_j}$) and so $p_j$
vanishes on $A_2 \cap X$. Therefore
each
$p_j$ vanishes on $A_2$, and so $A_1 \subseteq A_2$.
By reversing the roles of $A_1$ and $A_2$, we get that
$A_2 \subseteq A_1$.
\ep
\bd
\label{def:1.1} 
{\rm
Let $A$ be an algebraic set in $\C^n$.
Let us agree to say $A$ is {\it toral}
if $\T^n$ is determining for $A$ and 
say $A$ is {\it atoral} if $\T^n$ is not determining for any of the irreducible
components of $A$.
If $p \in \C[z_1,\ldots,z_n]$,
let us agree to say $p$ is {\it toral} (resp, {\it atoral}) if the
algebraic set $Z_p$ is.
}
\ed

Note that if $V$ is a variety, then $V$ is either toral or atoral
and if, in addition, $V$ is nonempty, then $V$ cannot be both toral and atoral.
Also, the empty set is both toral and atoral.

The definition of ``is determined by" immediately implies that
if $A_1$ is determining for $B$ and $A_2$ is determining for
$B$, then $A_1 \cup A_2$ is determining for $B$.
Therefore, finite unions of toral sets are toral.
That finite unions of atoral sets are atoral
is also true and follows immediately from the 
definition. Notice that it also follows
immediately from the definition of atoral,
that if $A$ is an atoral set, then each 
irreducible component of $A$ is atoral.
The corresponding assertion for toral 
algebraic sets is also true as the following
proposition shows.

\begin{prop}
\label{prop:toral:implies:components:toral}
If  $A$ is a toral algebraic set,
then each irreducible component of $A$ is toral.
\end{prop}
\new{This proposition follows from the following lemma.}

\begin{lemma}
\label{lemma:see:components}
If  $X$ is a set in $\C^n$, $A$ is an algebraic set in $\C^n$,
and $X$ is determining for $A$, 
then $X$ is determining for each irreducible component of $A$.
\end{lemma}
\bp
Suppose that $A$ is a toral algebraic set and $C$ is an irreducible
component of $A$. Set $A_0 = A \backslash C$.
\new{
Let $f \in Hol(C)$ and assume that $f |_{C \cap X} = 0$.}
There exists a function $\chi \in Hol(A)$ such that
$\chi |_{A_0} = 0$ and $\chi |_C \not\equiv 0$.
If $z \in A_0 \cap C$, then $\chi(z) f(z) = 0$.
Therefore, the function defined by
setting
$g(z) = \chi(z) f(z)$ for $z \in C$ and $g(z) = 0$ for $z \in A_0$
is well defined.

To show $g \in Hol(A)$, let $\zeta \in A$. We seek to show that there is a
neighborhood $U$ of $\zeta$ and an analytic function defined on
$U$ which equals $g$ on $U \cap A$.

Suppose $\zeta \in C\setminus \overline{A_0}$.
There exists an open disk $U$ centered at $\zeta$ such that
$U \cap A_0 = \empty$, and there exist
$\Gamma, F \in Hol(U)$ such that
$\Gamma |_{U \cap A} = \chi |_{U \cap A}$
and
$F |_{U \cap A} = f |_{U \cap A}$.
Therefore, $g |_{U \cap A} = \chi f |_{U \cap A} 
= \Gamma F |_{U \cap A}$.

Now suppose $\zeta \in A_0\backslash C$.
There exists an open disk $U$ centered at $\zeta$ such that
$U \cap C = \empty$. If $G(z) = 0$ for all $z \in U$, then
$g |_{U \cap A} = G |_{U \cap A}$.

Now suppose $\zeta \in \overline{A_0} \cap C$.
There exists an open disk $U$ centered at $\zeta$
and
$\Gamma, F \in Hol(U)$ such that
$\Gamma |_{U \cap A} = \chi |_{U \cap A}$
and
$F |_{U \cap C} = f |_{U \cap C}$.
Let $G(z) = F(z) \Gamma(z)$ for $z \in U$.
For $z \in A_0 \cap U$,
$G(z) = \Gamma(z) F(z) = \chi(z) f(z) = 0 = g(z)$ by
the choice of $\chi$.
For $z \in C \cap U$, $g(z) = f(z) \chi(z)$
and $G(z) = \Gamma(z) F(z) = \chi(z) f(z) = g(z)$.
Therefore, $g |_{U \cap A} = G |_{U \cap A}$.

We have shown that
$g \in Hol(A)$.
\new{
Note that $g |_{A \cap X } = 0$ and therefore, $g = 0$
since $X$ is determining for $A$ and $g \in Hol(A)$.}
Therefore, $\chi f |_C \equiv 0$.
Since $C$ is a variety, $Hol(C)$ is an integral domain.
Since $\chi f |_C = 0$ and $\chi |_C \not=0$,
$f |_C = 0$. Therefore, $C$ is toral.
\ep

If $B$ is an algebraic set,  
we define the {\it toral component of $B$} to be the union
of the irreducible toral components of $B$
and define
the {\it atoral component of $B$} to be the union
of the irreducible atoral components of $B$.

Let $\C^* = \C \backslash\{0\}$ and,
for $\zeta \in (\C^*)^n$, let
$\Flip{\zeta} = \left(\Flip{\zeta_1},\ldots,\Flip{\zeta_n}\right)$.
For $\zeta \in \C^n$ and
$d = (d_1,\ldots,d_n)$, an $n$-tuple of nonnegative integers,
let $\zeta^d = \zeta_1^{d_1} \zeta_2^{d_2}\cdots \zeta_n^{d_n}$.
Note that $\Flip{\zeta} = \zeta$  if and only if $\zeta \in \T^n$
and that $\zeta^d \in \C^*$ whenever $\zeta \in (\C^*)^n$.

Let us agree to say that an algebraic set $A \subseteq \C^n$
is \tnsym\ if,
for all $\zeta \in \cstarn$,
$\zeta \in A$  implies
$\Flip{\zeta} \in A$.
Note that,
since $\Flip{\Flip{\zeta}}=\zeta$,
$A$ is \tnsym\ if and only if,
for all $\zeta \in \cstarn$,
$\zeta \in A$  if and only if
$\Flip{\zeta} \in A$.

\bprop
\label{thm:toral:algebraicsets:are:symmetric}
If $A$ is a toral 
algebraic set in $\C^n$, then $A$ is \tnsym.
\eprop
\bp
Let $A$ be a toral
algebraic set in $\C^n$.
If $A = \C^n$, then $A$ is \tnsym.
Suppose $A \not= \C^n$.
There exist nonzero polynomials $p_1,\ldots,p_m$
such that $A = Z_{p_1} \cap Z_{p_2} \cap \cdots \cap Z_{p_m}$.
Let $q_j(z) = z^{\deg p_j} \overline{p_j(\Flip{z})}$.
Each $q_j$ is a nonzero polynomial.
For $\zeta \in A \cap \T^n$, 
$q_j(\zeta) = \zeta^{\deg p_j} \overline{p_j(\zeta)} = 0$.
Therefore, $q_j |_{A \cap \T^n} = 0$.
Since $A$ is toral,
$q_j | A \in Hol(A)$,
and $q_j |_{A \cap \T^n} = 0$,
we have
$q_j |_A = 0$.

If $\zeta \in A \cap (\C^*)^n$, then,
for each $j$,
$q_j(\zeta) = 0$ and $p_j(\Flip{\zeta}) = 0$.
Therefore, $\Flip{\zeta} \in A$.
Thus, $A$ is \tnsym.
\ep

In the special case when $n=2$, it is an easy matter
to describe the torality (resp, atorality)
of any algebraic set in terms of the 
torality (resp, atorality) of a single polynomial.
In $\C^2$, 
varieties are either points,
$\C^2$, or $Z_p$ for some irreducible $p \in \C[z_1,z_2]$, 
so any algebraic set $A$ in $\C^2$ can be represented in the
form $A = F \cup Z_p$ for some $p \in \C[z_1,z_2]$
where $F$ is the finite set of isolated       
points of $A$.
(If $A$ is finite, choose $p(z_1,z_2)=1$.)
The following proposition
follows from this characterization of irreducible
algebraic sets in $\C^2$.
\begin{prop}
\label{prop:alg:geom:two}
Let $A$ be an algebraic set in $\C^2$.
If $A$ is a finite set, then
$A$ is toral if and only if $A\subset \T^2$ and
$A$ is atoral if and only if $A\cap \T^2 = \empty$.
If $A$ is not a finite set and we let $A = F \cup Z_p$
where $F$ is the set of isolated
points of $A$ and $p \in \C[z_1,z_2]$, then  the following statements hold.
\begin{description}
\item{(a)}{ $A$ is toral if and only if
$p$ is toral and $F \subset \T^2$.}
\item{(b)}{ $A$ is atoral if and only if
$p$ is atoral and $F \cap\T^2 = \empty$.}
\end{description}
\end{prop}

\section{Toral and Atoral polynomials in $\C^n$ }
\label{secc2}

Let us agree to say that two irreducible polynomials $p$ and $q$     
are {\it essentially equal} if $p= c q$ for some nonzero $c \in \C$.
The zero set $Z_p$ is a variety if and only if $p$ is irreducible.
Since 
the zero set of a nonzero polynomial $p$ is equal to the union of the zero sets of
its irreducible factors, the results of Section \ref{secc} imply the following corollaries.

\begin{cor}
\label{cor:divisors:toral:are:toral}
Let $p$ be a nonzero polynomial in $\C[z_1,\ldots,z_n]$.
The following are equivalent.\\
\n \phantom{\ \ }i){ $p$ is toral (respectively, atoral)}\\
\n \phantom{\ }ii) {each irreducible factor of $p$ is toral (respectively, atoral)}\\
\n iii){ every divisor of $p$ is toral (respectively, atoral).}
\end{cor}

\begin{cor}
\label{cor:toral:atoral:factorization}
Let $p$ be a nonzero polynomial in $\C[z_1,\ldots,z_n]$.
There exist $q$, $r \in \C[z_1,\ldots,z_n]$ such that
$p=qr$, $q$ is toral and $r$ is atoral.
Moreover, if $p = q'r'$ is another factorization      
with $q'$ toral and $r'$ atoral, then
$q'$ is essentially equal to $q$
and
$r'$ is essentially equal to $r$.
\end{cor} 

\bp
The existence of $q$ and $r$ follows from
factoring $p$ into irreducible factors and then
grouping the toral and atoral factors.
Suppose that $p = q'r'$ is another factorization
with $q'$ toral and $r'$ atoral.
By Corollary \ref{cor:divisors:toral:are:toral},
every irreducible factor of $q'$ is toral
and every irreducible factor of $r'$ is atoral.
Thus, $q'$ divides $q$ and $r'$ divides $r$.
Since $qr=q'r'$, $q'$ is essentially equal to $q$
and
$r'$ is essentially equal to $r$.
\ep

Any nonzero polynomial in $\czn$ can be reflected in $\tn$ in the following way.
Let $p$ be a polynomial of degree $d=(d_1,\ldots,d_n)$.
We define the polynomial $\tp$ by
\be
\label{defa1}
\tp (z) \= z^d \ \overline{p(1/\overline{z})} \, .
\ee
Notice that $(pq)^\sim = p^\sim q^\sim$ and that 
$p^{\sim\sim} = p$ if and only if none of the
coordinate functions divide $p$.

We shall say that a polynomial $p$ is {\em $\T^n$-symmetric}
if $\tp = p$ or $p$ is the zero polynomial, and {\em essentially $\T^n$-symmetric} if there is a
unimodular constant $\tau$ such that $\tau p$ is $\T^n$-symmetric.

If $p$ is an irreducible atoral polynomial, 
then there is some non-zero $f$ in $Hol(Z_p)$ that
vanishes on $Z_p \cap \tn$. 
Therefore $Z_p \cap \tn$ is contained in the $(n-2)$-dimensional
analytic set $Z_p \cap Z_f$. In Theorem~\ref{thmb1} we show that
$Z_p \cap \tn$ is contained in an $(n-2)$-dimensional
\va\ (in other words, $f$ can be chosen to be a polynomial).
This can be thought of as an algebraic characterization of whether
a polynomial is toral: measuring directly how fat $Z_p \cap \tn$ is.

\bt
\label{thmb1}
If $p \in \czn$, then $p$ is 
atoral if and only if $Z_p \cap \tn$ is contained in an algebraic set $A$ of dimension
$n-2$.
\et
\bp
By Corollary
\ref{cor:divisors:toral:are:toral}, it is sufficient to consider polynomials
$p$ that are irreducible.

(Sufficiency)
Let
$$
A \= \bigcup_{j=1}^m \ \left( Z_{q_{1j}} \cap Z_{q_{2j}} \right) ,
$$
where $q_{1j}$ and $q_{2j}$
are relatively prime polynomials, and suppose $A$ contains
$Z_p \cap \tn$. Relabelling, if necessary, we can
assume that $p$ does not divide any $q_{1j}$.
If  $q = \Pi_j\ q_{1j}$, then $q$ vanishes on $Z_p \cap \tn$, but not on all $Z_p$.
Thus $p$ is atoral.
\bs
(Necessity)
If $p$ is not essentially $\T^n$-symmetric, then 
$$
Z_p \cap \tn \ \subseteq \ Z_p \cap Z_\tp,
$$
and we are done. So, assume that $p$ is $\T^n$-symmetric.
We shall show that the gradient of $p$ must vanish on 
$Z_p \cap \tn$.
Indeed, let $(z_0',w_0)$ be a point in 
$Z_p \cap \tn$ with $z_0' \in \C^{n-1}$, and assume that 
$$
\frac{\partial p}{\partial z_n} |_(z_0',w_0) \ \neq \ 0.
$$
By the implicit function theorem, there is an open neighborhood
$U$ of $z_0'$, an open neighborhood $W$ of $w_0$, 
and a holomorphic function $h$ on $U$ such that
if $(z',w)$ is in $U \times W$, then
\be
\label{eqc1}
p(z',w) = 0 \ \iff \ w = h(z') .
\ee
By (\ref{eqc1}) for every $z' \inn U \cap \T^{n-1}$, $p$ has only
one root $w$ in $W$. Since $p$ is $\T^n$-symmetric,
this root must be unimodular. Therefore,
$$
Z_p \cap \tn \cap ( U \times W) 
\= 
\{ (z',h(z') ) \ : \ z' \inn U \cap \T^{n-1} \} .
$$
Now suppose $f \inn Hol(Z_p)$ vanishes on $Z_p \cap \tn$.
The function $z' \mapsto f(z',h(z'))$ on $U$ vanishes
on $U \cap \T^{n-1}$, and therefore vanishes identically on $U$.
\new{Therefore, if we let $S$ be the set of singular points of $Z_p$, then 
$f$ vanishes on an open subset of $Z_p\backslash S$, and as $Z_p\backslash S$ is
connected
\cite{fis01},
$f$ is identically zero. }
Thus $p$ must be toral, a contradiction.

Therefore we can let $\dis A = Z_p \cap 
Z_{\frac{\partial p}{\partial z_n}}$.
\ep

\bs

\bprop
\label{propc1}
Every toral polynomial is essentially $\T^n$-symmetric.
\eprop
\bp
By Corollary 
\ref{cor:divisors:toral:are:toral}, it is sufficient to show that
every irreducible toral polynomial is $\T^n$-symmetric.
Since $Z_p \cap \tn = Z_{p^\sim}\cap \tn$, Theorem
\ref{thmb1} implies $p^\sim$ is toral.
Since $p^\sim$ vanished on $Z_p \cap \tn$ and $p$ is toral,
we must have $Z_p \subseteq Z_{p^\sim}$.
Since $p$ vanished on $Z_{p^\sim} \cap \tn$ and $p^\sim$ is toral,
we must have $Z_{p^\sim} \subseteq Z_p$.
Thus $Z_p = Z_{p^\sim}$ and since both $p$ and $p^\sim$ are irreducible,
we must have that one is a unimodular constant times the other.
\ep

\bs
Now we give a geometric condition which is sufficient to guarantee
that a polynomial is toral.
\bt
\label{thmb2}
Let $p \inn \czn$, and suppose $Z_p$ is disjoint from $\D^n \cup
\E^n$. Then $p$ is toral.
\et
\bp
%
Suppose $z' \inn \tnm$, and $p(z',w) = 0$.
Then $|w| =1$, or else some small perturbation of $(z',w)$ would 
yield a zero of $p$ in either $\dn$ or $\en$.

Write $$
p(z',w) \= \sum_{j=0}^k a_j (z') w^j .
$$
Let $D(z') \inn \cznm$ be the discriminant of 
$p(z',w)$; then $D$ will vanish
precisely at those points $z'$ such that the function
$w \mapsto p(z',w)$ has a root of multiplicity higher than one
for some $w$.

Let $B = Z_{D} \cup Z_{a_k}$. 
This is an algebraic set in $\cnm$ that does not disconnect $\cnm$.
Off $B$, one can choose $k$ holomorphic functions
$ \{ w_j (z') \}_{j=1}^k $
that take values in the $k$ sheets of $Z_p$ over $\cnm \setminus
B$.

Now let $f \inn Hol(Z_p)$ vanish on $Z_p \cap \tn$.
Locally, around any point $(z',w)$ in $Z_p$ with $z' \ \notin \ B$,
$f$ can be written as a function of $z'$, and this function
vanishes on $\tnm \setminus B$. By analytic continuation, $f$ must
vanish at any point $(z',w)$ with $z' \inn \C \setminus B$, and
by continuity of $f$, it must vanish on all of $Z_p$.
Therefore $Z_p$ is toral, as required.
\ep

\section{Inner functions}
\label{secd}

W.~Rudin showed that any
 rational inner function can be represented as 
$$
\phi(z)
\= z^h \frac{\tp(z)}{p(z)}
$$
for some polynomial $p$ that does not vanish on $\dn$, and
some monomial $z^h$, \cite[Thm. 5.2.5]{rud69}.
We show that the atoral factor of $p$ uniquely determines $\phi$.
Note that in the case $n=2$, B.~Cole and J.~Wermer have obtained
additional information about the relation between $p$ and $\tp$
\cite{colwer99}.
\bt
\label{thmd1}
Let 
\be
\phi(z) = z^h \frac{\tp(z)}{p(z)}
\qquad\mbox{\ and\  } \qquad
\psi(z) = z^e \frac{\tq(z)}{q(z)}
\ee
be two rational inner functions, with $p$ and $q$ polynomials that
do not vanish on $\dn$. Then $\phi$ and $\psi$ are essentially equal if
and only if
$p$ and $q$ 
have the same atoral factor and $h = e$.
\et
\bp
By Proposition~\ref{propc1}, any toral factor $r$ of either $p$ or $q$
is essentially $\T^n$-symmetric, so $r^\sim / r$ is constant. Thus we can
assume that both $p$ and $q$ are atoral. Moreover, if $p$ had any
nonconstant
$\T^n$-symmetric divisor $r$, then $Z_r$ would be disjoint from $\dn$ (since
$Z_p$ is) and from $\en$ (by $\T^n$-symmetry). Therefore by
Theorem~\ref{thmb2}, $r$ would be toral. So we can assume that neither
$p$ nor $q$ has any nonconstant $\T^n$-symmetric divisors.
 
To show that $p$ and $q$ must then be essentially equal, cross-multiply to
get a scalar $\tau$ such that 
\be
\label{eqd1}
\tau
z^h \tp(z) q(z) \= z^e \tq(z) p(z).
\ee
Since $p(0) \neq 0 \neq q(0)$, both $\tp$ and $\tq$ have 
$\overline{p(0)} z^{{\rm deg\ }p}$
and
$\overline{q(0)} z^{{\rm deg\ }q}$
respectively as their highest order terms.
Therefore the degree of the left-hand side of (\ref{eqd1}) is
$h + {\rm deg\ } p + {\rm deg\ }q$, and the 
degree of the right-hand side is
$e + {\rm deg\ } q + {\rm deg\ }p$. Therefore $h = e$.

Now since $p$ has no nonconstant $\T^n$-symmetric divisors, $p$ is relatively prime to
 $\tp$ (for if
$r$ were an irreducible polynomial that divided both, either $r$
would
be essentially $\T^n$-symmetric, or $r r^\sim$ would be a $\T^n$-symmetric 
polynomial that divided $p$). Therefore, $\tp$ divides $\tq$, and,
since $p=p^{\sim\sim}$ and $q=q^{\sim\sim}$,  $p$
divides $q$. Interchanging the roles of $p$ and $q$,
$p$ and $q$ must be essentially equal.
\ep
\bs

The following proposition shows that the zero set of a rational inner function is atoral, and the level set
for any unimodular number is toral.
Part (iii) is due to W.~Rudin \cite[Thm. 5.2.6]{rud69}.

\bprop
\label{propd2}
Let $\phi$ be a nonconstant rational inner function, and let $\alpha
\inn \C$. Then
\newline
\noindent
(i) If $\alpha \inn \D \cup \E$, then $Z_{\phi - \alpha}$ is atoral.
\newline
\noindent
(ii) If $\alpha \inn \T$, then $Z_{\phi - \alpha}$ is toral.
\newline
\noindent
(iii) If $\alpha \inn \D$, then $Z_{\phi-\alpha} \cap \E^n = \empty$.
\newline
\noindent
(iv) If $\alpha \inn \E$, then $Z_{\phi-\alpha} \cap \D^n = \empty$.
\newline
\noindent
(v) If $\alpha \inn \T$, then $Z_{\phi-\alpha} \cap \D^n = Z_{\phi-\alpha}  \cap \E^n = \empty$.
\eprop
\bp
(i,iii,iv) Write 
\be
\label{eqd3}
\phi(z) \=
 z^h \frac{\tp(z)}{p(z)},
\ee
where $p$ (and therefore $\tp$) is atoral.
Suppose first that $\alpha = 0$. Then $Z_\phi = Z_{z^h} \cup Z_{\tp}$.
The set $Z_{z^h}$ is disjoint from $\tn$, so is atoral, and $Z_{\tp}$
is atoral by the choice of $p$. Moreover, since $Z_p$ is disjoint from 
$\D^n$, $Z_{p^\sim}$ is disjoint from $\E^n$. Thus, 
$Z_{\phi} \cap \E^n = \empty$.

Now, if $\alpha \inn \D$, consider
$$
\psi(z) \= \frac{\phi(z) - \alpha}{1 - \overline{\alpha} \phi(z)} .
$$
Then $\psi$ is rational and inner, $Z_\psi = Z_{\phi - \alpha}$,
and $Z_\psi \cap \E^n = \empty$.

Finally, let $\alpha \inn \E$. 
Since $\phi$ is an inner function, the maximum principle,
$Z_{\phi-\alpha} \cap \D^n=\empty$. Also,
since $\alpha \inn \E$,
$\phi(\zeta) = \alpha$ if and only if $\phi(1/\overline{\zeta}) =
1/\overline{\alpha}$. So the zero set of $\phi - \alpha$ is the
reflection of the zero set of $\phi-1/\overline{\alpha}$, and therefore is
atoral.
\vs
(ii,v)
Suppose $|\alpha | = 1$. Off $Z_p$, 
we have $\phi(\zeta) = \alpha$ if and only if $\phi(1/\overline{\zeta}) =
1/\overline{\alpha} = \alpha$.
So $Z_{\phi - \alpha}$ is $\T^n$-symmetric. $Z_{\phi-\alpha}$
is disjoint from $\D^n$ by the
maximum principle, so $Z_{\phi-\alpha}$ must also be disjoint from $\en$. Therefore, by
Theorem~\ref{thmb2}, the set is toral.
\ep
\bs
The singular set $S_\phi$ of a rational inner function is the set
of points on $\T^n$ to which the function cannot be continuously
extended from $\D^n$.
If the function has the form
(\ref{eqd3}), it is the set $Z_p \cap \tn$.
If $\zeta$ is in this singular set, then $\tp(\zeta) = \zeta^{h+d}
\overline{p(\zeta)} = 0$, so 
$S_\phi \subseteq Z_p \cap Z_\tp$. 
Therefore we have:
\bprop
\label{propd3}
The singular set of a rational inner function is always
contained in an algebraic set of dimension $n-2$.
\eprop

\section{Application to Interpolation}
\label{secf}

Let $\hib$ denote the Banach algebra of bounded analytic functions
on the bidisk.
A {\em solvable Pick problem on $\D^2$} is a set $\{\l_1, \dots,
\l_N \}$ of points in $\D^2$ and a set $\{w_1, \dots, w_N\}$ of
complex numbers such that there is some function $\phi$ of norm
less
than or equal to one
in $\hib$ that interpolates (satisfies $\phi(\l_i) = w_i \ \forall
\ 1 \leq i \leq N$).
An {\em extremal Pick problem} is a solvable Pick problem for which
no function of norm less than one interpolates.
The points $ \l_i $ are called the {\em nodes},
and $w_i$ are called the {\em values}. By {\em interpolating
function} we mean any function in the closed
unit ball of $\hib$ that interpolates.

\bs
Consider the two following examples, in the case $N=2$.

Example 1. Let $\l_1 = (0,0),\, \l_2 = (1/2,0),\,
w_1 =0,\, w_2 = 1/2$. Then a moment's thought reveals that the
interpolating function is unique, and is given by
$\phi(z,w) = z$.

Example 2. Let $\l_1 = (0,0),\, \l_2 = (1/2,1/2),\,
w_1 =0,\, w_2 = 1/2$. Then the interpolating function is far from
unique --- either coordinate function will do, as will any convex
combination of them. (A complete description of all solutions is
given by J.~Ball and T.~Trent in \cite{baltre98}).
But on the algebraic set $\{(z,z) \, :
\,
z \, \in \, \D \}$, all solutions coincide by Schwarz's lemma.
\bs
For an arbitrary solvable Pick problem, let $\U$ be the set of
points in $\D^2$
on which all the interpolating functions in the closed unit ball of
$\hib$ have the same
value. The preceding examples show that $\U$ may be either the
whole bidisk or a proper subset. In the event that $\U$ is not the
whole bidisk, it is an algebraic set intersected with $\D^2$. 
Indeed, for any $\l_{N+1}$ not in
$\U$, there are two distinct values $w_{N+1}$ and $w'_{N+1}$
so that the corresponding $N+1$ point Pick problem has a solution.
By \cite{baltre98,agmc_bid} these problems have interpolating
functions that are rational, of degree bounded by $2(N+1)$. The
set $\U$ must lie in the zero set of the difference of these
rational functions. Taking the intersection over all $\l_{N+1}$ not
in $\U$, one gets that $\U$ is the intersection of the zero sets of
polynomials.  Since $\czn$ is Noetherian \cite{eis95},
$\U$ is the intersection of the zero sets of a finite
number of polynomials.
Therefore $\U$ is an algebraic set, and
indeed, by factoring these polynomials into their irreducible
factors, we see that $\U$ is the intersection with the bidisk of
the zero set of one polynomial, together with possibly a finite
number of isolated points. We shall call $\U$ the {\em uniqueness
set}. (If the problem is not extremal, $\U$ is just the
original set of nodes).

We shall say that an $N$-point extremal Pick problem is {\it
minimal} if
none of the $(N-1)$ point subproblems is extremal.
In \cite{agmc_dv}, a set $W$ was called a
 {\em distinguished variety} if $W$ the non-empty intersection of
the zero set of a polynomial with the bidisk, 
and moreover it satisfied the property
$$
\overline{W} \cap \partial (\D^2) \= \overline{W} \cap \T^2 .
$$
The following theorem was proved in \cite{agmc_dv}:
\bt
\label{thmf1}
The 
uniqueness variety of a minimal extremal Pick problem on $\D^2$
contains a distinguished variety that contains all the nodes.
\et
This theorem left open the possibility that $\U$ might still have
some isolated points in $\D^2$. We show that this cannot happen.
Indeed, $\U$ must be a toral algebraic set intersected with $\D^2$.
\bt
\label{thmf2}
The uniqueness set of a minimal extremal Pick problem on $\D^2$
has the form $\D^2 \cap Z_p$ where $p$ is a toral polynomial.
\et
We shall need to use Lojasiewicz's Vanishing theorem in the proof.
See \cite[Thm. 6.3.4]{kp02} for a proof of this form of the
theorem:
\bt
\label{thmf3} [Lojasiewicz] Let $f$ be a non-zero real analytic
function on an open set $U$ in $\R^n$. Assume that the zero set of
$f$ in $U$ is non-empty. Let $E$ be a compact subset of $U$. Then
there are constants $C$ and $k$ such that $|f(x)| \geq C\ {\rm
dist}(x,Z)^k$ for every $x \inn E$.
\et
{\sc Proof of Thm.~\ref{thmf2}:}
If $\U = \D^2$, take the polynomial to be $0$.
Otherwise,
let $\phi_1,\dots,\phi_m$ be rational inner functions that solve
the Pick problem and such that
$$
\U \= \left( \bigcap_{i \neq j} Z_{\phi_i - \phi_j} \right) \cap \D^2
.
$$
(This can be done because $\C[z_1,z_2]$ is Noetherian).
Let
$$
B \= \bigcap_{i \neq j} Z_{\phi_i - \phi_j} ,
$$
and let $V$ be the union of the irreducible toral components $V_i$ of $B$.
Notice that by Theorem~\ref{thmf1}, $V$ contains all the nodes of
the interpolation problem. 
Let
\be
\label{eqfa1}
\psi \= \frac{1}{m} ( \phi_1 + \dots + \phi_m) \= q/p 
\ee
be a rational solution, with $p$ and $q$ coprime polynomials,
and $q$ normalized to have modulus less than or equal to $1$ on
$\D^2$.
Let $$
S \= \bigcup_{j=1}^m S_{\phi_j} 
$$
be the union of the singular sets; then $S$ is finite by
Proposition~\ref{propd3}. 
As $p$ and $q$ are coprime, any zero of $p$ is a singularity of
$\psi$, so
$$
Z_p \cap \T^2 \ \subseteq \ S.
$$

Notice that $\psi$ will be unimodular on
$\T^2\setminus S$ 
only when all of the $\phi_j$'s are equal. 
Therefore, 
we have
$$
(B \cap \T^2) \, \cup \, S
 \= \{ \tau \inn \T^2 \ : \ |\psi(\tau)| = 1 \}  \, \cup \, S.
$$
For each $\tau \inn \T^2$, define
$$
l_\tau (z) \= 2 - \overline{\tau}_1 z_1 - \overline{\tau}_2 z_2 ,
$$
a linear polynomial whose only zero in $\overline{\D}^2$ is at
$\tau$. 
For each irreducible component $V_i$ of $V$, let $r_i$ be an
irreducible polynomial that vanishes on $V_i$. 
Define
\be
\label{eqf2}
g \= \left( \Pi\, r_i \right) \ \Pi \{ l_\tau \ : \
\tau \inn \left[ (B \setminus V) \cap \T^2 \right] \cup S
\} .
\ee
Since we are working in 2 dimensions, Theorem~\ref{thmb1} implies that atoral
varieties intersect $\T^2$ in a finite set. 
As $B \setminus V$ is contained in the atoral component of $B$,
we see that $( B
\setminus V ) \cap \T^2$ is finite. Thus, the
second product in (\ref{eqf2}) is over a finite set and $g$ is a
polynomial.
Furthermore, we have 
\be
\label{eqf3}
Z_g \cap \T^2 \= (B \cup S ) \cap \T^2 .
\ee
Now $|p|^2 - |q|^2$ is strictly greater than $0$ on
$\overline{\D^2} \setminus (B \cup S)$.
So applying Theorem~\ref{thmf3}
to the real analytic function $|p|^2 - |q|^2$ on $\T^2$,
this function must grow at least as fast as some power of the
distance to its zero set. Since $g|_{\T^2}$ vanishes 
on $(B \cup S) \cap \T^2$, and $g$ is a polynomial, we know that
$|g|^2$ can grow no faster than a constant times the distance
to $(B \cup S) \cap \T^2$.
Therefore
we conclude that there exist constants $\vare \, >
\, 0$ and $M \inn {\mathbb N}$ such that
\be
\label{eqf4}
2 \vare |g|^M + \vare^2 |g |^{2M}\  \leq \ |p|^2 - |q|^2
\quad {\rm on\ } \T^2 .
\ee
Now, let
$$
h \= \vare g^M .
$$
With this definition of $h$, we have that 
$h$ is zero on the nodes of the interpolation problem
(since $g$ vanishes on $V$, and $V$ contains these nodes by
Theorem~\ref{thmf1}),
$\| \psi + \psi h \| \leq 1$ (by (\ref{eqf4})), 
and
$h \neq 0$ on $\D^2 \setminus V$ (since by 
(\ref{eqf2}) the zeroes of $h$ that are not in $V$ are a union of
hyperplanes that just graze the closed bidisk at a single point).

Therefore $\psi + \psi h$ also solves the interpolation problem,
and so
\be
\label{eqf7}
\U \ \subseteq\ Z_{\psi h}   \cap \D^2 \=
\left( Z_\psi \cup V \right) \cap \D^2.
\ee

Suppose now that there is some point $\l$ in
$\left( Z_\psi \setminus V \right) \cap \D^2$.
This means that $\sum \phi_j(\l) = 0$, but not every $\phi_j(\l)$
is $0$. Replace $\psi$ in 
(\ref{eqfa1}) by $\psi_\l$,
some other strict convex combination of the $\phi_j$'s
that has the additional property that $\psi_\l (\l) \neq 0$.
Now repeat the above argument with $\psi_\l$ instead of $\psi$, and
in place of (\ref{eqf7}) we get
\be
\label{eqf8}
\U \ \subseteq\ Z_{\psi h}   \cap \D^2 \=
\left( Z_{\psi_\l} \cup V \right) \cap \D^2.
\ee
As
$$
\left[ \bigcap_{\l \in Z_\psi \setminus V} \
Z_{\psi_\l} \right] \cap 
\left[ Z_\psi \, \setminus V \right] \cap \D^2 
\= \emptyset ,
$$
combining (\ref{eqf7}) and (\ref{eqf8}) we get
$$
\U \ \subseteq\ V \cap \D^2 .
$$
As 
$$
B \cap \D^2 \ \subseteq \ \U,
$$
we conclude that 
$$
\U \ =\ V \cap \D^2,
$$
By Proposition \ref{prop:alg:geom:two}, $\U$ equals the intersection of the 
bidisk and the 
zero set of a toral polynomial,
as desired.

\ep

\bibliography{references}

\begin{thebibliography}{10}

\bibitem{agmc_dv}
J.~Agler and J.E. M\raise.45ex\hbox{c}Carthy.
\newblock Distinguished varieties.
\newblock {\em Acta Math.}
\newblock To appear.

\bibitem{agmc_bid}
J.~Agler and J.E. M\raise.45ex\hbox{c}Carthy.
\newblock {Nevanlinna-Pick} interpolation on the bidisk.
\newblock {\em J. Reine Angew. {M}ath.}, 506:191--204, 1999.

\bibitem{and63}
T.~{And\^o}.
\newblock On a pair of commutative contractions.
\newblock {\em Acta Sci. Math. (Szeged)}, 24:88--90, 1963.

\bibitem{bsv04}
J.A. Ball, C.~Sadosky, and V.~Vinnikov.
\newblock Conservative linear systems, unitary colligations and {Lax-Phillips}
  scattering: multidimensional generalizations.
\newblock {\em Internat. J. Control}, 77(9):802--811, 2004.

\bibitem{baltre98}
J.A. Ball and T.T. Trent.
\newblock Unitary colligations, reproducing kernel {Hilbert} spaces, and
  {Nevanlinna-Pick} interpolation in several variables.
\newblock {\em J. Funct. Anal.}, 197:1--61, 1998.

\bibitem{bv00}
J.A. Ball and V.~Vinnikov.
\newblock Hardy spaces on a finite bordered {R}iemann surface, multivariable
  operator theory and {F}ourier analysis along a unimodular curve.
\newblock In {\em Operator Theory Advances and Applications}, volume 129, pages
  37--56. {Birkh\"auser}, Basel, 2000.

\bibitem{bv03}
J.A. Ball and V.~Vinnikov.
\newblock Overdetermined multidimensional systems: state space and frequency
  domain methods.
\newblock In {\em Mathematical systems theory in biology, communications,
  computation, and finance}, volume 134 of {\em IMA Vol. Math. Appl.}, pages
  63--119. Springer, Berlin, 2003.

\bibitem{colwer99}
B.J. Cole and J.~Wermer.
\newblock {And\^o's} theorem and sums of squares.
\newblock {\em Indiana Math. J.}, 48:767--791, 1999.

\bibitem{eis95}
D.~Eisenbud.
\newblock {\em Commutative Algebra}.
\newblock Springer, New York, 1995.

\bibitem{fis01}
G.~Fischer.
\newblock {\em Plane algebraic curves}.
\newblock American Mathematical Society, Providence, 2001.

\bibitem{kp02}
S.~Krantz and H.~Parks.
\newblock {\em A primer of real analytic functions}.
\newblock {Birkh\"auser}, Basel, 2002.

\bibitem{rud69}
W.~Rudin.
\newblock {\em Function Theory in {Polydiscs}}.
\newblock Benjamin, New York, 1969.

\bibitem{szn-foi}
B.~Szokefalvi-Nagy and C.~Foia\c{s}.
\newblock {\em Harmonic Analysis of Operators on {Hilbert} Space}.
\newblock North Holland, Amsterdam, 1970.

\end{thebibliography}

\end{document}